\documentclass[a4paper,12pt]{article}
\usepackage{amsfonts,amssymb,latexsym,amsmath}
\usepackage{multicol}
\usepackage[cp1251]{inputenc}
\usepackage[russian]{babel}
\begin{document}
\renewcommand{\refname}{References}

\newcommand{\Dp}{\Delta^{\!+}}
\newcommand{\De}{\Delta}

\newcommand{\gog}{{\mathfrak g}}
\newcommand{\nog}{{\mathfrak n}}
\newcommand{\mog}{{\mathfrak m}}
\newcommand{\hog}{{\mathfrak h}}
\newcommand{\bog}{{\mathfrak b}}
\newcommand{\pog}{{\mathfrak p}}
\newcommand{\rog}{{\mathfrak r}}
\newcommand{\ut}{{\mathfrak u}{\mathfrak t}}
\newcommand{\gl}{{\mathfrak g}{\mathfrak l}}

\newcommand{\sg}{\sigma}
\newcommand{\eps}{\varepsilon}
\newcommand{\ad}{{\mathrm{ad}}}
\newcommand{\Ad}{{\mathrm{Ad}}}
\newcommand{\id}{{\mathrm{id}}}
\newcommand{\GL}{\mathrm{GL}}
\newcommand{\UT}{\mathrm{UT}}
\newcommand{\ord}{\mathrm{ord}}

\newcommand{\AC}{{\cal A}}
\newcommand{\AG}{{\frak A}}
\newcommand{\BC}{{\cal B}}
\newcommand{\SC}{{\cal S}}
\newcommand{\UC}{{\cal U}}
\newcommand{\FC}{{\cal F}}
\newcommand{\YC}{{\cal Y}}
\newcommand{\Xb}{{\Bbb X}}
\newcommand{\Nb}{{\Bbb N}}

\date{}
\title{Regular $N$-orbits in the nilradical of\\ a parabolic subalgebra}
\author{A. N. Panov\thanks{This research was partially supported by the RFBR (projects 05-01-00313, 06-01-00037)}
\and V. V. Sevostyanova}

\maketitle

\begin{center}
\parbox[b]{330pt}{\small\textsc{Abstract.} In the present paper the
adjoint action of the unitriangular group in the nilradical of a
parabolic subalgebra is studied. We set up general conjectures on
the construction of the field of invariants and the structure of
orbits of maximal dimension. The conjecture is proved for parabolic
subalgebras of special types.}
\end{center}

\vspace{0.5cm}

Consider the general linear group $\GL(n,K)$ defined over an
algebraically closed field $K$ of characteristic 0. Let $B$ ($N$,
respectively) be its Borel (maximal unipotent, respectively)
subgroup, which consists of triangular matrices with nonzero (unit,
respectively) elements on the diagonal. We fix a parabolic subgroup
$P$ that contains $B$. Denote by $\pog$, $\bog$ and $\nog$ the Lie
subalgebras in $\gl(n,K)$ that correspond to $P$, $B$ and $N$
respectively. We represent $\pog=\rog\oplus\mog$ as the direct sum
of the nilradical $\mog$ and a block diagonal subalgebra $\rog$ with
sizes of blocks $(n_1,\ldots, n_s)$. The subalgebra $\mog$ is
invariant relative to the adjoint action of the group $P$,
therefore, $\mog$ is invariant relative to the action of the
subgroups $B$ and $N$. We extend this action to the representation
in the algebra $\AC=K[\mog]$ and in the field $\FC=K(\mog)$. The
subalgebra $\mog$ contains a Zariski-open $P$-orbit, which is called
the Richardson orbit. Consequently, the algebra of invariants
$\AC^P$ coincides with $K$. Invariants of the adjoint action of the
group $N$ in $\mog$ are studied worse. In the case $P=B$, the
algebra of invariants $\AC^N$ is the polynomial algebra
$K[x_{12},x_{23},\ldots,x_{n-1,n}]$. Let $\rog$ be the sum of two
blocks; this case is a result of~\cite{B}. We do not know when the
algebra of invariants $\AC^N$ is finitely generated.

We formulate a number of conjectures (Conjectures 1--3) on the
structure of the field of invariants $\FC^N=\mathrm{Fract}\AC^N$ and
on the descriptions of regular $N$-orbits (i.e. $N$-orbits of
maximal dimension). We prove Conjectures 1 and 2 for parabolic
subalgebras of special types (Theorems 1--2). Theorem 3 provides a
description of the field $\FC^B$ for the same parabolic subalgebras.
Propositions 1 and 2 give a partial decision of Conjecture 3. We
construct the system of generators of the algebra $\AC^N$ for a
parabolic subalgebra with sizes of blocks $(2,4,2)$ (Proposition 3).

We begin with definitions. Every positive root $\gamma$ in
$\gl(n,K)$ has the form~(see \cite{GG}) $\gamma=\eps_i-\eps_j$,
$1\leqslant i<j\leqslant n$. We identity a root $\gamma$ with the
pair $(i,j)$ and the set of positive roots $\Dp$ with the set of
pairs $(i,j)$, $i<j$. The system of positive roots $\Dp_\rog$ of the
reductive subalgebra $\rog$ is a subsystem in $\Dp$.

Let $\{E_{ij}:~i<j\}$ be the standard basis in $\nog$. By $E_\gamma$
denote the basis element $E_{ij}$, where $\gamma=(i,j)$.

We define a relation in $\Dp$ such that $\gamma'\succ\gamma$
whenever $\gamma'-\gamma\in\Dp_\rog$. If $\gamma\prec\gamma'$ or
$\gamma\succ\gamma'$, then the roots $\gamma$ and $\gamma'$ are
comparable. Denote by $M$ the set of $\gamma\in\Dp$ such that
$E_\gamma\in\mog$. We identify the algebra $\AC$ with the polynomial
algebra in the variables $x_{ij}$, ~$(i,j)\in M$.

\medskip
\textbf{Definition 1.} A subset $S$ in $M$ is called a \emph{base}
if the elements in $S$ are not pairwise comparable and for any
$\gamma\in M\setminus S$ there exists $\xi\in S$ such that
$\gamma\succ\xi$.

\medskip
Note that $M$ has a unique base $S$, which can be constructed in the
following way. We form the set $S_1$ of minimal elements in $M$ (we
say that $\gamma$ is a minimal element in $S_1$ if there is no
$\xi\in S_1$ such that $\gamma\succ\xi$). By definition, $S_1\subset
S$. We form a set $M_1$, which is obtained from $M$ by deleting
$S_1$ and all elements
$$\{\gamma\in M:\exists\ \xi\in S_1,\ \gamma\succ\xi\}.$$
The set of minimal elements $S_2$ in $M_1$ is also contained in $S$,
and so on. Continuing the process, we get the base $S$.

\medskip
\textbf{Definition 2.} An ordered set of positive roots
$\{\gamma_1,\ldots,\gamma_s\}$ is called a \emph{chain} if
$\gamma_1=(a_1,a_2)$, $\gamma_2=(a_2,a_3)$,
$\gamma_3=(a_3,a_4)$,\ldots

\medskip
\textbf{Definition 3.} We say that two roots $\xi,\xi'\in S$ form an
\emph{admissible pair} $q=(\xi,\xi')$ if there exists
$\alpha_q\in\Dp_\rog$ such that the ordered set of roots
$\{\xi,\alpha_q,\xi'\}$ is a chain. Note that the root $\alpha_q$ is
uniquely determined by $q$.

\medskip
We form the set $Q:=Q(\pog)$ that consists of admissible pairs of
roots in $S$. For every admissible pair $q=(\xi,\xi')$ we construct
a positive root $\varphi_q=\alpha_q+\xi'$. Consider the subset
$\Phi=\{\varphi_q:~ q\in Q\}$.

Let $\mathfrak{p}$ be any parabolic subalgebra. We construct a
diagram by $\mathfrak{p}$, which is a square $n\times n$-matrix.
Roots from $S$ are marked by symbol $\otimes$ and roots from $\Phi$
are labeled by the symbol $\times$ in the diagram. The other entries
in the diagram are empty.

Let a parabolic subalgebra $\mathfrak{p}$ be the subalgebra of type
$(2,1,3,2)$ (the type of a parabolic subalgebra is the sizes of
diagonal blocks). We have the following diagram.
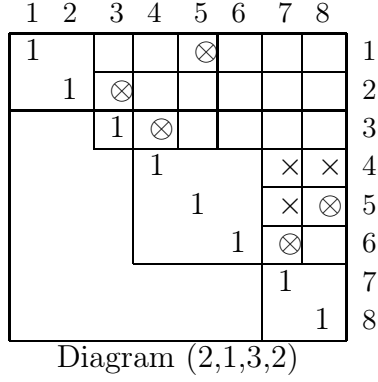
\begin{center}
\begin{tabular}{|p{0.1cm}|p{0.1cm}|p{0.1cm}|p{0.1cm}|p{0.1cm}|
p{0.1cm}|p{0.1cm}|p{0.1cm}|l} \multicolumn{2}{l}{{\small
1\hspace{5pt} 2\ }}&\multicolumn{2}{l}{{\small 3\hspace{5pt}
4}}&\multicolumn{2}{l}{{\small 5\hspace{5pt} 6}}&
\multicolumn{2}{l}{{\small 7\hspace{5pt} 8}}\\
\cline{1-8} \multicolumn{2}{|l|}{1}&&&$\otimes$&&&&{\small1}\\
\cline{3-8} \multicolumn{2}{|r|}{1}&$\otimes$&&&&&&{\small2}\\
\cline{1-8} \multicolumn{2}{|c|}{}&1&$\otimes$&&&&&{\small3}\\
\cline{3-8} \multicolumn{3}{|c|}{}&\multicolumn{3}{|l|}{1}&$\times$&$\times$&{\small4}\\
\cline{7-8} \multicolumn{3}{|c|}{}&\multicolumn{3}{|c|}{1}&$\times$ &$\otimes$&{\small5}\\
\cline{7-8} \multicolumn{3}{|c|}{}&\multicolumn{3}{|r|}{1}&$\otimes$&&{\small6}\\
\cline{4-8} \multicolumn{6}{|c|}{}&\multicolumn{2}{|l|}{1}&{\small7}\\
\multicolumn{6}{|c|}{}&\multicolumn{2}{|r|}{1}&{\small8}\\
\cline{1-8} \multicolumn{8}{c}{Diagram (2,1,3,2)}\\
\end{tabular}
\end{center}

Consider the formal matrix $\Xb$ in which the variables $x_{ij}$
occupy the positions $(i,j)\in M$ and the other entries are equal to
zero. For any root $\gamma=(a,b)\in M$ we denote by $S_\gamma$ the
set of $\xi=(i,j)\in S$ such that $i>a$ and $j<b$. Let
$S_\gamma=\{(i_1,j_1),\ldots,(i_k,j_k)\}$. Denote by $M_\gamma$ a
minor $M_I^J$ of the matrix $\Xb$ with the ordered systems of rows
$I=\ord\{a,i_1,\ldots,i_k\}$ and columns $J=\ord\{j_1,\ldots,j_k,
b\}$.

For every admissible pair $q=(\xi,\xi')$, we construct the
polynomial
\begin{equation}
L_q=\sum_{\scriptstyle\alpha_1,\alpha_2\in\Dp_\rog\cup\{0\}
\atop\scriptstyle\alpha_1+\alpha_2=\alpha_q}
M_{\xi+\alpha_1}M_{\alpha_2+\xi'}.\label{L_q}
\end{equation}

\medskip
\textbf{Conjecture 1.} The field of invariants $\FC^N$ is the field
of rational functions of polynomials $M_\xi$, $\xi\in S$, and $L_q$,
$q\in Q$.

\medskip
A next conjecture is a consequence of the preceding one.

\medskip
\textbf{Conjecture 2.} The maximal dimension of an $N$-orbit in
$\mog$ is equal to $\dim\mog-|S|-|Q|$.

\medskip
Denote by $\YC:=\YC_\pog$ the subset in $\mathfrak{m}$ that consists
of matrices
$$\sum_{\xi\in S}c_{\xi}E_\xi+\sum_{\varphi\in\Phi}c'_{\varphi}E_\varphi.$$

\medskip
\textbf{Conjecture 3.} Any regular $N$-orbit (i.e., an orbit of
maximal dimension) has a nonzero intersection with $\YC$.

\medskip
\textbf{Notation.} We can replace the set $\Phi$ by any similar
subset $\Psi$ in the following way. We can replace the root
$\alpha_q+\xi'$ in $\Phi$ by one of two roots $\xi+\alpha_q$ and
$\alpha_q+\xi'$.

\medskip
\textbf{Theorem 1.} For an arbitrary parabolic subalgebra, the
system of polynomials
$$\{M_\xi,~\xi\in S,~L_q,~q\in Q,\}$$
is contained in $\AC^N$ and is algebraically independent over
$K$.

\medskip
\textsc{Proof.} The representation of $P$ in $\AC=K[\mog]$ is
determined by $T_gf(x) = f(\Ad_g^{-1}x),$ where $g\in P$ and
$f\in\AC$. The action of $T_g$ in $\AC$ is uniquely defined by the
action on $x_{i,j}$, $(i,j)\in M$. The elements $x_{i,j}$ make the
matrix $T_g\Xb=g^{-1}\Xb g$, where the formal matrix $\Xb$ is
defined above.

A polynomial $f$ of $\AC$ is an $N$-invariant if $f$ is an invariant
of the adjoint action of any one-parameter subgroup
$g_k(t)=1+tE_{k,k+1}$, $1\leqslant k <n$. The action of $g_i(t)$ on
the matrix $\Xb$ reduces to the composition of two transformations:
\begin{itemize}
\item[1)] the row with number $k+1$ multiplied by $-t$ is added to
the row with number $k$,
\item[2)] the column with number $k$ multiplied by $t$ is added to
the column with number $k+1$.
\end{itemize}

The invariance of $M_{(a,b)}$ follows from the notations:
\begin{itemize}
\item[1).] Numbers of rows and columns of the minor $M_{(a,b)}$ fill
the segments of natural numbers $I=[a,\max I]$, $J=[\min J,b]$.
\item[2).] All elements $(i,j)$ of the matrix $\Xb$ are equal to zero,
where $a\leqslant i\leqslant n$ and $1\leqslant j<\min J$ or $\max
I< i\leqslant n$ and $1\leqslant j\leqslant b$.
\end{itemize}

Now let us prove that $L_q$ is in $\AC^N$. This statement follows
from the invariance of $L_q$ under the adjoint action of the
one-parameter subgroup $g_k(t)$, where $g_k(t)$ corresponds to the
simple root $\beta=(k,k+1)$, $1\leqslant k< n$.

Let $q=(\xi,\xi')$, where $\xi=(a,b)$,~$\xi'=(a',b')$. Using the
definition of admissible pair, we have $a<b<a'<b'$ and
$\alpha_q=(b,a')\in\Dp_\rog$. If $k<b$ or $ k\geqslant a'$, then the
minors of the right part of (\ref{L_q}) are $g_k(t)$-invariants.

If $b\leqslant k<k+1\leqslant a'$, then
$\alpha_q=\gamma_1+\beta+\gamma_2$, where
$\gamma_1,\gamma_2\in\Dp_\rog\cup\{0\}$. We have
\begin{equation}
\left\{\begin{array}{l}
T_{g_k(t)}M_{\xi+\gamma_1+\beta}=
M_{\xi+\gamma_1+\beta}+tM_{\xi+\gamma_1},\\
T_{g_k(t)}M_{\beta+\gamma_2+\xi'}=
T_{\beta+\gamma_2+\xi'}-tM_{\gamma_2+\xi'}.\\
\end{array}\right.\label{T_g(M_xi)}
\end{equation}
The other minors of (\ref{L_q}) are invariants under the action of
$g_k(t)$. Combining (\ref{L_q}) and (\ref{T_g(M_xi)}), we get
$$\left(T_{g_k(t)}L_q\right)-L_q=M_{\xi+\gamma_1}
\left(M_{\beta+\gamma_2+\xi'}-t M_{\gamma_2+\xi'}\right)+$$
$$\left(M_{\xi+\gamma_1+\beta}+t
M_{\xi+\gamma_1}\right)M_{\gamma_2+\xi'}-
M_{\xi+\gamma_1}M_{\beta+\gamma_2+\xi'}-
M_{\xi+\gamma_1+\beta}M_{\gamma_2+\xi'}=0.$$

To prove the second statement of the theorem, we need an order
relation on the set of roots $S\cup\Phi$ such that
\begin{itemize}
\item[1)] $\xi<\varphi$ for any $\xi\in S$ and $\varphi\in\Phi$;
\item[2)] for other pairs of roots from $S\cup\Phi$, the relation $<$ means
the lexicographic order relation.
\end{itemize}
Consider the restriction homomorphism $\pi:f\mapsto f|_\YC$ of the
algebra $\AC$ to $\YC$. The image of $\AC$ is the polynomial algebra
$K[\YC]$ of $x_\xi$, $\xi\in S$, and of $x_\varphi$,
$\varphi\in\Phi$. The image of $M_\xi$, $\xi\in S$, has the form
\begin{equation}
\pi\left(M_\xi\right)=\displaystyle\pm x_\xi\prod_{\xi'\in
S_{\xi}}x_{\xi'}.\label{pi(M)}
\end{equation}
Suppose the root $\varphi\in\Phi$ corresponds to the admissible pair
$q=(\alpha,\beta)$. Then the image of the polynomial $L_q$ has the
form
\begin{equation}
\pi\left(L_{q}\right)=\displaystyle\pm
x_{\varphi}x_{\alpha}\prod_{\alpha'\in
S_{\alpha}}x_{\alpha'}\prod_{\beta'\in
S_{\alpha}}x_{\alpha'}.\label{pi(L)}
\end{equation}

The system of the images $\left\{\pi\left(M_\xi\right), \xi\in S,
\pi\left(L_q\right), q\in Q\right\}$ is algebraically independent
over $K$. Therefore, the system $\left\{M_\xi, \xi\in S, L_q, q\in Q
\right\}$ is algebraically independent over $K$.~$\Box$

\medskip
Consider the open subset $\UC_0=\{x\in\mog:~M_\xi\ne
0,~\forall\xi\in S\}\subset\mog$.

\medskip
\textbf{Proposition 1.} Let $\mathfrak{p}$ be a parabolic subalgebra
of type $\left(n_1,\ldots,n_s\right)$. Suppose
$n_1\geqslant\ldots\geqslant n_s$; then the $N$-orbit of any
$x\in\mathcal{U}_0$ intersects $\mathcal{Y}$ at a unique point.

\medskip
\textsc{Proof.} Let $n_1\geqslant\ldots\geqslant n_s$, then the
number of elements of $S$ is equal to $n_2+\ldots+n_s$. The set $S$
consists of the roots $\xi_{i,j}=(m_i-j+1,m_i+j)$, where
$m_i=n_1+\ldots+n_i$ and $1\leqslant i\leqslant s-1$, $1\leqslant
j\leqslant n_{i+1}$.

\begin{enumerate}
\item
Let us show that for any $A=\left(a_{ij}\right)\in\mathcal{U}_0$
there exists $g\in N$ such that $\mathrm{Ad}_gA\in\mathcal{Y}$. The
proof is by induction on the number of diagonal blocks. Suppose that
the statement is true for a parabolic subalgebra $\pog_1$, where
$\pog_1\subset\gl(n-n_1,K)$ has the reductive subalgebra of type
$(n_2\geqslant\ldots\geqslant n_s)$. Let us show that the statement
is true for a parabolic subalgebra $\pog$ of type $(n_1\geqslant
n_2\geqslant\ldots\geqslant n_s)$.

Consider the Lie algebra $\gl(n-n_1,K)$ regarded as a subalgebra of
$\gl(n,K)$. Let $\gl(n-n_1,K)$ have zeros in the first $n_1$ rows
and columns. The systems $S_1\subset S$ and $\Phi_1\subset\Phi$
correspond to the parabolic subalgebra $\pog_1$. For the algebra
$\pog$ of type $(2,2,2,1,1)$ and the subalgebra $\pog_1$ of type
$(2,2,1,1)$, we have the following diagrams.
\begin{center}
\begin{tabular}{|p{0.1cm}|p{0.1cm}|p{0.1cm}|p{0.1cm}|p{0.1cm}|
p{0.1cm}|p{0.1cm}|p{0.1cm}|l} \multicolumn{2}{l}{{\small
1\hspace{5pt} 2\ }}&\multicolumn{2}{l}{{\small 3\hspace{5pt}
4}}&\multicolumn{2}{l}{{\small 5\hspace{5pt} 6}}&
\multicolumn{2}{l}{{\small 7\hspace{5pt} 8}}\\
\cline{1-8} \multicolumn{2}{|l|}{1}&&$\otimes$&&&&&{\small1}\\
\cline{3-8} \multicolumn{2}{|r|}{1}&$\otimes$&&&&&&{\small2}\\
\cline{1-8} \multicolumn{2}{|c|}{}&\multicolumn{2}{|l|}{1}&$\times$&$\otimes$&&&{\small3}\\
\cline{5-8} \multicolumn{2}{|c|}{}&\multicolumn{2}{|r|}{1}&$\otimes$&&&&{\small4}\\
\cline{3-8} \multicolumn{4}{|c|}{}&\multicolumn{2}{|l|}{1}&$\times$&&{\small5}\\
\cline{7-8} \multicolumn{4}{|c|}{}&\multicolumn{2}{|r|}{1}&$\otimes$&&{\small6}\\
\cline{5-8} \multicolumn{6}{|c|}{}&1&$\otimes$&{\small7}\\
\cline{7-8} \multicolumn{7}{|c|}{}&1&{\small8}\\
\cline{1-8} \multicolumn{8}{c}{Diagram (2,2,2,1,1)}\\
\end{tabular}
\hspace{1.5cm}
\begin{tabular}{|p{0.1cm}|p{0.1cm}|p{0.1cm}|p{0.1cm}|p{0.1cm}|
p{0.1cm}|p{0.1cm}|p{0.1cm}|l} \multicolumn{2}{l}{{\small
1\hspace{5pt} 2\ }}&\multicolumn{2}{l}{{\small 3\hspace{5pt}
4}}&\multicolumn{2}{l}{{\small 5\hspace{5pt} 6}}&
\multicolumn{2}{l}{{\small 7\hspace{5pt} 8}}\\
\cline{1-8} \multicolumn{2}{|c|}{}&\multicolumn{6}{|c|}{}&{\small1}\\
\multicolumn{2}{|c|}{}&\multicolumn{6}{|c|}{}&{\small2}\\
\cline{1-8} \multicolumn{2}{|c|}{}&\multicolumn{2}{|l|}{1}&&$\otimes$&&&{\small3}\\
\cline{5-8} \multicolumn{2}{|c|}{}&\multicolumn{2}{|r|}{1}&$\otimes$&&&&{\small4}\\
\cline{3-8} \multicolumn{4}{|c|}{}&\multicolumn{2}{|l|}{1}&$\times$&&{\small5}\\
\cline{7-8} \multicolumn{4}{|c|}{}&\multicolumn{2}{|r|}{1}&$\otimes$&&{\small6}\\
\cline{5-8} \multicolumn{6}{|c|}{}&1&$\otimes$&{\small7}\\
\cline{7-8} \multicolumn{7}{|c|}{}&1&{\small8}\\
\cline{1-8} \multicolumn{8}{c}{Diagram (2,2,1,1)}\\
\end{tabular}
\end{center}

By the inductive assumption, all elements $a_{i,j}$ of a matrix $A$
from the nilradical of $\pog$ are equal to zero for all
$n_1<i,j\leqslant n$ and $(i,j)\not\in S_1\cup\Phi_1$.

Consider the matrix $\mathrm{Ad}_{g_1}A$, where
$$g_1=\mathrm{exp}\left(1+
t_1E_{n_1+1,n_1+2}+\ldots+t_{n_{2}-1}E_{n_1+1,n_1+n_2}+\right.$$
$$\left.+t'_1E_{n_1-1,n_1}+\ldots+t'_{n_1-1}E_{1,n_1}\right).$$
By the condition, $a_{n_1,n_1+1}\ne 0$. There exist
$t_1,\ldots,t_1',\ldots$ such that the elements of the matrix
$\mathrm{Ad}_{g_1}A$ at the entries $(1,n_1+1),\ldots,(n_1-1,n_1+1)$
and $(n_1,n_1+2),\ldots (n_1,n_1+n_2)$ are equal to zero.

In general, elements of the matrix $\mathrm{Ad}_{g_1}A$ at the
entries $(i,n_1+n_2+1)$ are not equal to zero, where $n_1+1\leqslant
i\leqslant n_1+n_2-1$. The entries $(i,n_1+n_2+1)$ fill by the
symbol $\times$ (see the entry (3,5) in the Diagram $(2,2,2,1,1)$).

Similarly, we find $g_2,\ldots,g_{n_2}$ such that all elements of
the matrix
$$A'=\mathrm{Ad}_{g_{n_2}}\ldots
\mathrm{Ad}_{g_{2}}\mathrm{Ad}_{g_{1}} A$$ are equal to zero, except
for the elements of $S\cup\Phi$ or for the block
\begin{equation}
\{(i,j):~1\leqslant i\leqslant n_1,\ n_1+n_2+1\leqslant j\leqslant
n\}.\label{block}
\end{equation} The
elements from the block (\ref{block}) are labeled by the symbol~$*$
in the following diagram.
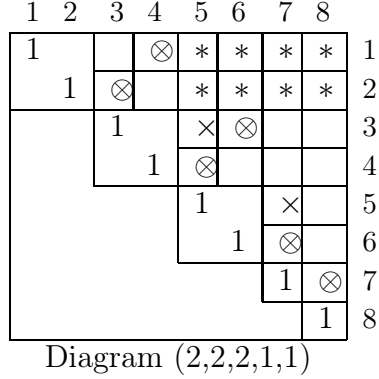
\begin{center}
\begin{tabular}{|p{0.1cm}|p{0.1cm}|p{0.1cm}|p{0.1cm}|p{0.1cm}|
p{0.1cm}|p{0.1cm}|p{0.1cm}|l} \multicolumn{2}{l}{{\small
1\hspace{5pt} 2\ }}&\multicolumn{2}{l}{{\small 3\hspace{5pt}
4}}&\multicolumn{2}{l}{{\small 5\hspace{5pt} 6}}&
\multicolumn{2}{l}{{\small 7\hspace{5pt} 8}}\\
\cline{1-8} \multicolumn{2}{|l|}{1}&&$\otimes$&$*$&$*$&$*$&$*$&{\small1}\\
\cline{3-8} \multicolumn{2}{|r|}{1}&$\otimes$&&$*$&$*$&$*$&$*$&{\small2}\\
\cline{1-8} \multicolumn{2}{|c|}{}&\multicolumn{2}{|l|}{1}&$\times$&$\otimes$&&&{\small3}\\
\cline{5-8} \multicolumn{2}{|c|}{}&\multicolumn{2}{|r|}{1}&$\otimes$&&&&{\small4}\\
\cline{3-8} \multicolumn{4}{|c|}{}&\multicolumn{2}{|l|}{1}&$\times$&&{\small5}\\
\cline{7-8} \multicolumn{4}{|c|}{}&\multicolumn{2}{|r|}{1}&$\otimes$&&{\small6}\\
\cline{5-8} \multicolumn{6}{|c|}{}&1&$\otimes$&{\small7}\\
\cline{7-8} \multicolumn{7}{|c|}{}&1&{\small8}\\
\cline{1-8} \multicolumn{8}{c}{Diagram (2,2,2,1,1)}\\
\end{tabular}
\end{center}

Note that the marked by the symbol $\otimes$ elements of the matrix
$A'$ are not equal to zero.

We shall show that there is the element $h\in N$ such that the
marked by~$*$ elements of the matrix $\mathrm{Ad}_hA'$ are equal to
zero. We start with the $n$th column. Let the symbol $\otimes$ be at
the entry $(i,n)$ of the last column for some~$i$. There exist
$s_1,\ldots,s_{n_1}$ such that the entries $(1,n),\ldots,(n_1,n)$ of
the matrix $\mathrm{Ad}_{h_1}A'$ are equal to zero, where
$$h_1=\mathrm{exp}\left(1+
s_1E_{1,i}+\ldots+s_{n_1}E_{n_1,i}\right).$$ In the same way, we
find $h_2,\ldots,h_{n-n_1}$ such that
$$\mathrm{Ad}_{h_{n-n_1}}\ldots\mathrm{Ad}_{h_2}A'\in\mathcal{Y}.$$

\item
Taking into account (\ref{pi(M)}) and (\ref{pi(L)}), we have that
the $N$-orbit of $A$ intersects $\mathcal{Y}$ at a unique
point.~$\Box$
\end{enumerate}

\medskip
\textbf{Proposition 2.} Let $\mathfrak{p}$ be a parabolic subalgebra
of type $\left(n_1,n_2,n_3\right)$, where $n_1,n_2,n_3$ are any
numbers. Then the $N$-orbit of any $x\in\mathcal{U}_0$ intersects
$\mathcal{Y}$ at a unique point.

\medskip
\textsc{Proof.} The proof is similarly.~$\Box$

\medskip
Let $\SC$ be the set of denominators generated by the minors
$M_\xi$, $\xi\in S$. We form the localization $\AC^N_\SC$ of the
algebra of invariants $\AC^N$ on $\SC$. Since the minors $M_\xi$ are
$N$-invariants, we have $\AC^N_\SC=\left(\AC_S\right)^N$.

\medskip
\textbf{Theorem 2.} Under the conditions of Proposition 1 or 2, we
have the following statements.
\begin{enumerate}
\item The ring $\AC^N_S$ is the ring of polynomials in
$M_\xi^{\pm 1}$, $\xi\in S$, and in $L_q$, $q\in Q$.
\item The field of invariants $\FC^N$ is the field of the rational
functions of $M_\xi$, $\xi\in S$, and $L_q$, $q\in Q$.
\end{enumerate}

\medskip
\textsc{Proof.} Consider the restriction homomorphism $\pi:f\mapsto
f|_\YC$ of the algebra $\AC^N$ to $K[\YC]$. Under the conditions of
Proposition 1 or 2, the image $\pi(M_\xi)$ is equal to the product
$$\pm x_\xi x_{\xi_1}\ldots x_{\xi_s},$$
where every $\xi_{i+1}$ is the greatest root in $S$, in the sense of
the above order, lesser than $\xi_i$. We extend $\pi$ to a
homomorphism
$$\pi_S:\AC^N_\SC\to K[\YC]_S,$$
where $K[\YC]_S$ is the localization of $K[\YC]$ with respect to
$x_\xi$, $\xi\in S$. We show that $\pi_S$ is an isomorphism.

If $f\in\mathrm{Ker}\,\pi_S$, then $f(\Ad_N\YC)=0$. Since, by
Proposition 1 and 2, $ \Ad_N\YC$ contains a Zariski-open subset,
then $f=0$. Consequently, $\pi_S$ is an embedding $\AC^N_\SC$ in
$K[\YC]_S$. Next, we write the formulates (\ref{pi(M)}) and
(\ref{pi(L)}) in the form
\begin{equation}
\begin{array}{l}
\pi\left(M_\xi\right)=\pm x_\xi\pi\left(M_{\xi_1}\right),\\
\pi\left( L_q\right)=\pm
x_\varphi\pi\left(M_{\xi}\right)\pi\left(M_{\xi'_1}\right),
\end{array}\label{local1}
\end{equation}
where the admissible pair $q=(\xi,\xi')$ corresponds to the root
$\varphi\in\Phi$ and $\xi_1$ ($\xi'_1$, respectively) is the
greatest root in the base that is less than $\xi$ ($\xi'$,
respectively) in the sense of the lexicographical order. We get
\begin{equation}
\begin{array}{l}
\pi_S\left(M_\xi M^{-1}_{\xi_1}\right)=\pm x_\xi,\\
\pi_S\left(L_q M^{-1}_{\xi} M^{-1}_{\xi'_1}\right)=\pm x_\varphi.\\
\end{array}\label{local2}
\end{equation}
From (\ref{local2}) it follows that the image of $\pi_S$ coincides
with $K[\YC]_S$. Thus, $\pi_S$ is an isomorphism.~$\Box$

\medskip
By $\AG$ denote the system of weights $\alpha_q$, ~$q\in Q$.

\medskip
\textbf{Theorem 3.} Under the conditions of Proposition 1 or 2, we
have the field $\FC^B$ is the field of the rational functions of
$\mathrm{corank}(\AG)$ variables.

\medskip
\textsc{Proof.} Consider the Cartan subgroup $H\subset\GL(n,K)$ of
the diagonal matrices. The field of invariants $\FC^B$ is a subfield
of $\FC^N$ and coincides with $\left(\FC^N\right)^H$. The system of
roots $S\cup\Phi$ generates the lattice of weights of the field
$\FC^N$. The field $\FC^B$ is a transcendental extension of $K$ and
$$\mathrm{tr\,deg}\,\FC^B=\mathrm{corank}(S\cup\Phi).$$

The system $S$ is linearly independent and
$$\mathrm{rank}(S\cup\Phi)=
\mathrm{rank}(S\cup\AG)=\mathrm{rank}(S)+\mathrm{rank}(\AG).$$
Further, $\mathrm{corank}(S \cup \Phi)=\mathrm{corank}(\AG).$~$\Box$

\medskip
Consider the parabolic subalgebra $\pog\subset\gl(8,K)$ with sizes
of blocks $(2,4,2)$. We give the complete description of the algebra
of invariants $\AC^N$. The base $S$ consists of the roots
$\alpha_1=(2,3)$, $\alpha_2=(1,4)$, $\beta_1=(6,7)$,
$\beta_2=(5,8)$. Any pair $(\alpha_i,\beta_j)$ is an admissible one.
The parabolic subalgebra $\pog$ corresponds to the following
diagram.
\begin{center}
\begin{tabular}{|p{0.1cm}|p{0.1cm}|p{0.1cm}|p{0.1cm}|p{0.1cm}|
p{0.1cm}|p{0.1cm}|p{0.1cm}|l} \multicolumn{2}{l}{{\small
1\hspace{8pt} 2\ }}&\multicolumn{2}{l}{{\small 3\hspace{8pt}
4}}&\multicolumn{2}{l}{{\small 5\hspace{8pt} 6}}&
\multicolumn{2}{l}{{\small 7\hspace{8pt} 8}}\\
\cline{1-8} \multicolumn{2}{|l|}{1}&&$\otimes$&&&&&{\small1}\\
\cline{3-8} \multicolumn{2}{|r|}{1}&$\otimes$&&&&&&{\small2}\\
\cline{1-8} \multicolumn{2}{|c|}{}&\multicolumn{4}{|l|}{1}&$\times$&$\times$&{\small3}\\
\cline{7-8} \multicolumn{2}{|c|}{}&\multicolumn{4}{|l|}{\hspace{18pt}1}&$\times$&$\times$&{\small4}\\
\cline{7-8} \multicolumn{2}{|c|}{}&\multicolumn{4}{|l|}{\hspace{35pt}1}& &$\otimes$&{\small5}\\
\cline{7-8} \multicolumn{2}{|c|}{}&\multicolumn{4}{|r|}{1}&$\otimes$&&{\small6}\\
\cline{3-8} \multicolumn{6}{|c|}{}&\multicolumn{2}{|l|}{1}&{\small7}\\
\multicolumn{6}{|c|}{}&\multicolumn{2}{|r|}{1}&{\small8}\\
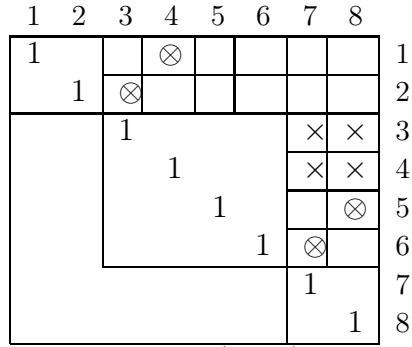
\cline{1-8} \multicolumn{8}{c}{Diagram (2,4,2)}\\
\end{tabular}
\end{center}

Let $M_1,M_2,N_1,N_2$ be the minors $M_{\alpha_1}$, $M_{\alpha_2}$,
$M_{\beta_1}$, $M_{\beta_2}$, respectively. By $L_{i,j}$ denote the
corresponding to the admissible pair $(\alpha_i,\beta_j)$
polynomial. By Theorem 2, the field of invariants $\FC^N$ is the
field of the the rational functions of $M_1$, $M_2$, $N_1$, $N_2$,
$L_{1,1}$, $L_{1,2}$, $L_{2,1}$, $L_{2,2}$. The generates have the
form
$$M_1=x_{23},\quad M_2=
\left\vert\begin{array}{cc}x_{13}&x_{14}\\x_{23}&x_{24}
\end{array}\right\vert,\quad N_1=x_{67},\quad N_2=
\left\vert\begin{array}{cc}x_{57}&x_{58}\\x_{67}&x_{68}
\end{array}\right\vert,$$
$$L_{1,1}=x_{23}x_{37}+x_{24}x_{47}+x_{25}x_{57}+x_{26}x_{67},$$
$$L_{1,2}=\left\vert\begin{array}{cc}x_{13}&x_{14}\\x_{23}&x_{24}
\end{array}\right\vert x_{47}
+\left\vert\begin{array}{cc}x_{13}&x_{15}\\x_{23}&x_{25}
\end{array}\right\vert x_{57} + \left\vert\begin{array}{cc}x_{13}&x_{16}\\x_{23}&x_{26}
\end{array}\right\vert x_{67},$$
$$L_{2,1}=x_{23}
\left\vert\begin{array}{cc}x_{37}&x_{38}\\x_{67}&x_{68}
\end{array}\right\vert + x_{24}
\left\vert\begin{array}{cc}x_{47}&x_{48}\\x_{67}&x_{68}
\end{array}\right\vert
+ x_{25}\left\vert\begin{array}{cc}x_{57}&x_{58}\\x_{67}&x_{68}
\end{array}\right\vert,$$
$$L_{2,2}=\left\vert\begin{array}{cc}x_{13}&x_{14}\\x_{23}&x_{24}
\end{array}\right\vert \cdot \left\vert\begin{array}{cc}x_{47}&x_{48}\\x_{67}&x_{68}
\end{array}\right\vert + \left\vert\begin{array}{cc}x_{13}&x_{15}\\x_{23}&x_{25}
\end{array}\right\vert \cdot \left\vert\begin{array}{cc}x_{57}&x_{58}\\x_{67}&x_{68}
\end{array}\right\vert.$$

By $D$ denote the minor $M^{7,8}_{1,2}$ of the matrix $\Xb^2$. It is
easily shown that $D$ is an $N$-invariant. We have the identity
\begin{equation}
L_{1,2}L_{2,1}-L_{1,1}L_{2,2}=M_1N_1D.\label{LL-LL=MMD}
\end{equation}

By $\BC_0$ denote the subalgebra such that the polynomials $M_i$,
$N_i$, $L_{i,j}$, $i,j=1,2$, generate $\BC_0$. By $\BC_1$ denote the
subalgebra such that $\BC_1$ is generated by $\BC_0$ and~$D$. Since
all generators are $N$-invariants, we have $\BC_0\subset\BC_1\subset
\AC^N$.

\medskip
\textbf{Proposition 3.} We have
\begin{itemize}
\item[1)] $\BC_0\ne\BC_1$;
\item[2)] $\AC^N=\BC_1$.
\end{itemize}

\medskip
\textsc{Proof.}
\begin{enumerate}
\item Suppose $\BC_0=\BC_1$, then the invariant $D$ is contained
in $\BC_0$. Therefore, there exists the polynomial
$f(u_1,\ldots,u_8)$ such that
\begin{equation}
D=f(M_1,M_2,N_1,N_2,L_{1,1},L_{1,2},L_{2,1},L_{2,2}).\label{D=f()}
\end{equation}
Combining (\ref{LL-LL=MMD}) and (\ref{D=f()}), we obtain that the
system $M_i$, $N_i$, $L_{ij}$, $i,j\in\{1,2\}$, is algebraically
dependent. This contradicts Theorem 1.

\item Let $\SC$ be the set of denominators generated by minors
$M_1$, $M_2$, $N_1$ and $N_2$. By Theorem 2, it follows that the
localization $\AC^N_\SC$ of the $N$-algebra of invariants on $\SC$
coincides with the algebra of Laurent polynomials
$$K[M_1^{\pm 1},M^{\pm 1}_2,N_1^{\pm 1},N_2^{\pm1},L_{11},L_{12},L_{21},
L_{22}].$$ If $f\in \AC^N$, then there exist $k_1,k_2,k_3,k_4$ such
that
$$M_1^{k_1}M_2^{k_2}N_1^{k_3}N_2^{k_4} f\in\BC_0.$$
Let us show that for any $M\in\{M_1,M_2,N_1,N_2\}$ we have that if
$F\in\AC$ and $MF\in\BC_1$, then $F\in\BC_1$. From this it follows
that $f\in\BC_1$.

We proof the theorem when $M=M_1$. The cases $M=M_2$, $M=N_1$ and
$M=N_2$ are similar. Let $F\in\AC^N$ and $M_1F\in\BC_1$. Denote
$M_1F=h$. We have $h|_{\mathrm{Ann}M_1}=0$. We form the matrix
$$Y:=Y_{a,b,c}:=
\left(\begin{array}{cccccccc}
0&0&a_1&0&0&0&0&0\\
0&0&0&a_2&0&0&0&0\\
0&0&0&0&0&0&c_{11}&c_{12}\\
0&0&0&0&0&0&c_{21}&c_{22}\\
0&0&0&0&0&0&0&b_2\\
0&0&0&0&0&0&b_1&0\\
0&0&0&0&0&0&0&0\\
0&0&0&0&0&0&0&0\\
\end{array}\right),$$
where $a_i$, $b_i$, $c_{ij}$ are any numbers, $i,j\in\{1,2\}$.

We have
$$M_1(Y)=0,\quad M_2(Y)=-a_1a_2,\quad N_1(Y)=b_1,\quad N_2(Y)=b_1b_2,$$
$$L_{11}(Y)=a_2c_{21},~L_{12}(Y)=-a_2b_1c_{22},$$$$
L_{21}(Y)=-a_1a_2c_{21},~L_{22}(Y)=a_1a_2b_1c_{22},$$
$$D(Y)=a_1a_2\left\vert
\begin{array}{cc}c_{11}&c_{12}\\
c_{21}&c_{22}\end{array}\right\vert=
a_1a_2(c_{11}c_{22}-c_{12}c_{21}).$$ By $V$ denote the space $K^9$
and by $K[V]$ denote the polynomial algebra
$K[u_1,u_2,v_1,v_2,w_{11},w_{12},w_{21}, w_{22},z]$. Since
$h\in\BC_1$, there exists the polynomial
$p(u_1,u_2,v_1,v_2,w_{11},w_{12},w_{21}, w_{22},z)\in K[V]$ such
that
\begin{equation}
h=p(M_1,M_2,N_1,N_2,L_{11},L_{12},L_{21},
L_{22},D).\label{h=p(M,N,L)}
\end{equation}
Since $h$ is equal to zero in $\mathrm{Ann}\,M_1$, then $h(Y)=0$. We
compute (\ref{h=p(M,N,L)}) at the point $Y$. We have
$$p(0,-a_1a_2,b_1,-b_1b_2,a_2c_{21},-a_2b_1c_{22},-a_1a_2c_{21},
a_1a_2b_1c_{22},$$$$a_1a_2(c_{11}c_{22}-c_{12}c_{21}))=0$$ for any
$a_i,b_j,c_{ij}\in K$. There exist $p_1,p_2\in K[V]$ such that
$$p=u_1p_1+(w_{12}w_{21}-w_{11}w_{22})p_2.$$ Combining the last
equation, (\ref{LL-LL=MMD}) and (\ref{h=p(M,N,L)}), we get
$$M_1F=M_1p_1+(L_{1,2}L_{2,1}-L_{1,2}L_{2,2})p_2=M_1p_1+M_1N_1Dp_2.$$
Hence, $F=p_1+N_1Dp_2\in\BC_1$.~$\Box$
\end{enumerate}

\textsc{Department of Mechanics and Mathematics, Samara State
University, Russia}\\ \emph{E-mail address}: \verb"apanov@list.ru",
\verb"victoria.sevostyanova@gmail.com"

\end{document}